\documentclass[10pt,twoside,reqno]{article}
\usepackage{amsmath}
\newtheorem{thm}{Theorem}[section]
\newtheorem{cor}[thm]{Corollary}
\newtheorem{lem}[thm]{Lemma}

\newtheorem{defn}[thm]{Definition}
\newtheorem{exa}[thm]{Example}

\numberwithin{equation}{section}

\begin{document}
\setcounter{page}{1}
\begin{center}
\vspace{0.4cm} {\large{\bf Fuzzy Rate Analysis of Operators  and its Applications in Linear Spaces}}\\
\vspace{0.4cm}
\vspace{0.4cm}

Yijin Zhang$^1$, Honggang Li$^{2*}$, Maoming Jin$^3$, Zongbing Lin$^4$\\

\vspace{0.4cm}

 $1$ Key Lab of Intelligent Analysis and Decision on Complex Systems,\\
   Key Laboratory of Industrial Internet of Things and  Networked Control, \\Ministry of Education,
Chongqing University of Posts and Telecommunications, Chongqing  400065, P.R. China.  E-mail: zhangyj@cqupt.edu.cn

$2$ Pass College of Chongqing Technology and Business University\\
 Chongqing 401520, P.R. China.  E-mail: 2273354482@qq.com  \\

$3$ Institute of Nonlinear Analysis Research,
Changjiang Normal University\\
 Chongqing 400803, P.R. China. E-mail: 19860025@yznu.cn \\

$4$ School of Mathematics and Computer Science, Panzhihua University,\\
 Panzhihua 617000, P.R. China. E-mail: zongbinglin@sohu.com \\

\vspace{0.4cm}

$*$ The corresponding author: Honggang  Li\\

\end{center}
\vspace{0.4cm}

\begin{abstract}  In this paper, a new concept, the fuzzy rate of an operator in linear spaces is proposed for the very first time.
Some properties and basic principles of it are studied. Fuzzy rate of an operator $B$ which is specific in a plane is discussed.
As its application, a new fixed point existence theorem is proved.


\end{abstract}

\section{Introduction}

More and more classical  analysis theory are being developed into fuzzy analysis theory. Fuzzy sets, fuzzy logic, fuzzy numbers,
fuzzy topologies and so on were introduced and studied[1-3]. Chang and Huang, Ding and Jong, Jin, Li and others studied several kinds of variational inequalities (inclusions) for
fuzzy mappings[4-8].

Recently, Konwar and Nabanita introduce the notion of continuous linear operators  and establish the uniform continuity theorem and Banach's contraction principle in an intuitionistic fuzzy n-normed linear space\cite{Konwar23}.Wang investigates the concepts and some properties of interval-valued fuzzy ideals in B-algebras and the homomorphic inverse image of interval-valued intuitionistic fuzzy ideals\cite{Wang24}. Fixed Point Theorems in Partially Ordered Fuzzy Metric Spaces and Operator Theory and Fixed Points in Fuzzy Normed Algebras and Applications are studied in
\cite{Cho}. Fuzzy-wavelet-like operators via a real-valued scaling function are discussed in \cite{Ezzati}. A linear fuzzy operator inequality approach is  proposed for the first time in \cite{Tai}. Fuzziness degree's quantity measure as to fuzzy operator is researched by means of  fuzzy set theory in \cite{Liu}. For more details, we
reference to the readers [1-15].

 In this work, we come up with  the concept of  fuzzy rate of an operator and consider its properties and applications. We also  explore fuzzy rate which is  produced by an operator effecting an element, as well as some  properties and applications of it.These are new extension, attempt and applications to  the operator theory in linear space and fuzzy theory.

The remainder of this paper is organized as follows. In Section 2, we  give an example which helps us introduce the concept of  fuzzy rate of an operator.  In Section 3, we propose
the concept and prove  some
basic properties of it.  In Section 4, a new Fixed Point Existence Theorem with the fuzzy
rate of the operator $B$ is obtained as its application.

\section{An Example}

 Here, an example is given to introduce a new concept, the fuzzy rate of an operator in linear spaces.

\begin{exa}
Let $U=R\times R$ be a real plane(Universe), $c$ be a cycle or ellipse whose center is at $(0, 0)$ on $U$, $F:U \to [0, 1]$  be a membership function, and

$\cal F$$(U)=\{F_c(x, y)|F_c(x, y)$ is the membership function for point $(x, y)$  belonging to a curve $c\}$.

In the plane $U$, we suppose that the equation of $c$ is
\begin{equation}
c_{(\lambda, r)}: x^2+\lambda y^2=r^2,
\end{equation}
where $r> 0$ and $0<\lambda $ are two parameters. Define a membership function $F _{c_{(1, r)}}(x, y)$   for  $(x, y)$  belonging to the curve $c_{(\mu, r)}(\mu> 0)$,

 \begin{equation} \label{e3.47}
F_{c_{(\mu, r)}}(x, y)=\begin{cases}
0, & (x, y)\notin c_{(\lambda, r)},\\
       e^{-(\lambda-\mu)^2}, & (x, y)\in c_{(\lambda, r)}(0<\lambda\neq \mu),\\
1, & (x, y)\in c_{(\lambda, r)}(0<\lambda= \mu),
 \end{cases}
\end{equation}

In the (2.2), let $\mu=1$, then the curve $c_{(1, r)}$ is the circle. We have $F_{c_{(1, r)}}(0, r)=1$ for the point $(0, r)\in c_{(1, r)}$ and $F_{c_{(1, r)}}(\frac{r}{\sqrt{2}}, \frac{r}{\sqrt{6}})=e^{-(3-1)^2}=e^{-4}$ for $\lambda=3$, but $F_{c_{(1, r)}}(r, r)=0$ for $(r, r)\notin c_{(\lambda, r)}$ and any $\lambda> 0$.

 Set
\begin{eqnarray}
B=\left(
                  \begin{array}{ccc}
                    1 &0\\
                    0 &b\\
                   \end{array}
                \right)
\end{eqnarray}
be an operator on $U$ for $b\neq 0$, then

$$F_{c_{(1, r)}}(B(0, r))=F_{c_{(1, r)}}(0, br)=e^{-( \frac{1}{b^2}-1)^2}.$$
where the point $B(0, r)\in c_{(\frac{1}{b^2}, r)}$ with the operator $B: (0, r)\mapsto B(0, r)$.

The  value,
\begin{eqnarray*}
\frac{F_{c_{(1, r)}}(B(0, r))}{F_{c_{(1, r)}}(0, r)}=\frac{e^{-( \frac{1}{b^2}-1)^2}}{1}=e^{-( \frac{1}{b^2}-1)^2}
\end{eqnarray*}

expresses a fuzzy rate of operator $B$ at the point $(0, r)\in c_{(1, r)}$ with $F_{c_{(1, r)}}$ for $\lambda=\frac{1}{b^2}.$

 Other examples,
if $r=1$, $b=\sqrt{2}$, then

\begin{eqnarray*}
\frac{F_{c_{(1, 1)}}(B(0, 1))}{F_{c_{(1, 1)}}(0, 1)}=\frac{e^{-( \frac{1}{b^2}-1)^2}}{1}=e^{-( \frac{1}{2}-1)^2}=e^{-\frac{1}{4}}.
\end{eqnarray*}
\end{exa}

 Furthermore, suppose that $\cal X$ is a linear space, $ {\O}\neq T \subseteq \cal X,$ $x \in \cal X$, $F_T(x):\cal X$$ \to [0, 1]$   is  a  membership function for  $x$ belonging to the set $T$[4, 5, 7, 8],
$B: \cal X \to \cal X $  is an operator. Then, $F_T(B(x))$ reflects the membership degree of the image of $x$  belonging to the set $T$. It's clear that the value $\frac{F_{T}(B(x))}{F_{T}(x)}$, the ratio of $F_T(B(x))$ and $F_T(x)$, indicates the changing rate which is produced by the mapping $B: x \to B(x) $. we can consider
a special value $\sup_{F\in \cal F (X)}\frac{F(B(y))}{F(y)}$ to express a fuzzy rate of operator $B$ at a point $y\in \cal X$ with $\cal F$$(\cal X)$,
it is  very interesting to consider the impact and properties of the operator with respect to $F_T$.

\section{Fuzzy Rate Analysis of Operators}

In this section, we first give the concept of a fuzzy rate of an operator,  then we show  some basic properties of it.

\begin{defn}
Let $\cal X$ be a linear space, $B:\cal X \to \cal X $  be  an operator, $\cal {B}( X)$$=\{B|B: \cal X\to \cal X\}$, $F:\cal X$ $ \to [0, 1]$ be a membership function over $\cal X$, $\cal P(X)$$=\{F|F:\cal X$$ \to [0, 1]\}$  be a collection of all membership functions over $\cal X$. For any $\emptyset\neq \cal F(X)\subseteq \cal P(X)$, if
\begin{equation}
{\|B\|_y}= \sup_{F\in \cal F (X)}\frac{ {F}(B(y))}{ {F}(y)}
\end{equation}
exists, then $ {\|B\|_y}$ is called a fuzzy rate of  the operator $B$ at the  point $y\in \cal X$ on $\cal F(X)$.
\end{defn}

For Example 2.1,let $\cal F$$(U)=\{F_{c_{(\mu, r)}}|\mu, r> 0\}$ by (2.2), we can achieve, obviously

\begin{eqnarray*}
&&{\|B\|}_{(0, r)}= \sup_{ {F_{c_{(\mu, r)}}\in {\cal F}(U)}}\frac{ {F_{c_{(\mu, r)}}(B(0, r))}}{ {F_{c_{(\mu, r)}}((0, r))}}\\
&&=\sup_{{F_{c_{(\mu, r)}}\in {\cal F}(U)}}\frac {F_{c_{(\mu, r)}}((0, br))} {F_{c_{(\mu, r)}}((0, r))}\\
&&=\sup_{{F_{c_{(\mu, r)}}\in {\cal F}(U)}}\frac {e^{-( \frac{1}{b^2}-\mu)^2}} { e^{-(1-\mu)^2}}\\
&&=\sup_{\mu> 0}e^{(1- \frac{1}{b^4})} e^{2\mu( \frac{1}{b^2}-1)}\\
&&=\begin{cases}
+\infty, & |b|< 1,\\
 
e^{(1- \frac{1}{b^4})}, & |b|\geq 1.
 \end{cases}
\end{eqnarray*}

At the same time, we have the following theorem about the relationship between a fuzzy rate of  the operator and fuzzy sets.

\begin{thm}
Let $\cal X$ be a linear space, $B:\cal X \to \cal X $  be  an operator, $\cal {B}( X)$$=\{B|B: \cal X\to \cal X\}$, $F:\cal X$ $ \to [0, 1]$ be a membership function over $\cal X$, $\cal P(X)$$=\{F|F:\cal X$$ \to [0, 1]\}$  be a collection of all membership functions over $\cal X$ and $ {\|B\|_y}$ be the  fuzzy rate of  the operator $B$ at the  point $y\in \cal X$ on $\cal F(X)$. Then for any $\emptyset\neq \cal F(X)\subseteq \cal P(X)$,
there exist two membership functions $ {F},  {G}\in \cal F (X)$ such that, for each  $y\in \cal X$,
\begin{equation}
 {\|B\|_y} {F}(y)=  {G}(B(y)).
\end{equation}
\end{thm}
\par

\textbf{Proof.} For any $\emptyset\neq \cal F (X)\subseteq \cal P(X)$, since
\begin{equation*}
 {\|B\|_y}= \sup_{ {F}\in \cal F(X)}\frac{ {F}(B(y))}{ {F}(y)}< +\infty,
\end{equation*}
then $ {F}(y)\neq 0$. $\forall n$, $\exists $ $ {G}_n\in \cal F(X)$, we arrive at
 \begin{equation*}
 {\|B\|_y}\geq \frac{ {G}_n(B(y))}{ {G}_n(y)}>  {\|B\|_y}-\frac{1}{n},
\end{equation*}
and
\begin{equation*}
 {\|B\|_y}=\lim_{n\to \infty} \frac{ {G}_n(B(y))}{ {G}_n(y)}< +\infty
\end{equation*}
for each $y\in \cal X$.

 Since $0<  {G}_n(y),  {G}_n(B(y))\leq 1$, there exist $ {G}_{n_k}(y)\rightarrow  {F}(y)$ as $k\rightarrow +\infty$, and $ {G}_{{n_k}_m}(B(y))\rightarrow  {G}(B(y))$ as $m\rightarrow +\infty$ for any $y\in \cal X$.

If $0<  {F}(y)\leq 1$, it follows that
\begin{equation*}
 {\|B\|_y}=\lim_{n\to \infty} \frac{ {G}_n(B(y))}{ {G}_n(y)}=\lim_{m\to \infty} \frac{ {G}_{{n_k}_m}(B(y))}{ {G}_{{n_k}_m}(y)}= \frac{\lim_{m\to \infty} {G}_{{n_k}_m}(B(y))}{\lim_{m\to \infty} {G}_{{n_k}_m}(y)}=\frac{ {G}(B(y))}{ {F}(y)},
\end{equation*}
and $ {\|B\|_y} {F}(y)=  {G}(B(y))$ for every $y\in \cal X$.

If $ {F}(y)=0$, it implies that
\begin{equation*}
 {\|B\|_y}=\lim_{n\to \infty} \frac{ {G}_n(B(y))}{ {G}_n(y)}=\lim_{m\to \infty} \frac{ {G}_{{n_k}_m}(B(y))}{ {G}_{{n_k}_m}(y)}= \frac{\lim_{m\to \infty} {G}_{{n_k}_m}(B(y))}{\lim_{m\to \infty} {G}_{{n_k}_m}(y)}<+\infty,
\end{equation*}
and $\lim_{m\to \infty} {G}_{{n_k}_m}(B(y))=0$ ,  $ {G}(B(y))=\lim_{m\to \infty} {G}_{{n_k}_m}(B(y))=0$\cite{Zorich25}.

Therefore, $ {\|B\|_y} {F}(y)=  {G}(B(y))$ holds for any $y\in \cal X$.

It is easy to  verify that  the converse proposition of Theorem 3.2 holds. We reach

\begin{thm}
Let $\cal X$ be a linear space, $B:\cal X \to \cal X $  be  an operator, $\cal {B}( X)$$=\{B|B: \cal X\to \cal X\}$, $F:\cal X$ $ \to [0, 1]$ be a membership function over $\cal X$, $\cal P(X)$$=\{F|F:\cal X$$ \to [0, 1]\}$  be a collection of all membership functions over $\cal X$. For any $\emptyset\neq \cal F (X)\subseteq \cal P(X)$, if there exist $ {F}$, $ {G}\in \cal F(X)$ and $ {F}(y)\neq 0$ such that
\begin{equation*}
 {\|B\|_y} {F}(y)=  {G}(B(y)),
\end{equation*}
then for any $y\in \cal X$,
\begin{equation*}
 {\|B\|_y}= \sup_{ {F}\in \cal F(X)}\frac{ {F}(B(y))}{ {F}(y)}=\frac{ {G}(B(y))}{ {F}(y)}< +\infty,
\end{equation*}
\end{thm}
that is, the  fuzzy rate of  the operator $B$ exists.

Now, we state some basic properties of the  fuzzy rate of  the operator $B$ as the next theorem. These  properties are very useful for further applications.

\begin{thm}
Let $\cal X$ be a linear space, $B:\cal X \to \cal X $  be  an operator, $\cal {B}( X)$$=\{B|B: \cal X\to \cal X\}$, $F:\cal X$ $ \to [0, 1]$ be a membership function over $\cal X$, $\cal P(X)$$=\{F|F:\cal X$$ \to [0, 1]\}$  be a collection of all membership functions over $\cal X$. For any $\emptyset\neq \cal F (X)\subseteq \cal P(X)$, $B_1, B_2\in \cal B(X)$ and  the  identity operator $I\in \cal B(X)$, then

(1) $ {\|B\|_y}> 0$ for any $B\in \cal B(X)$;

(2) $ {\|I\|_y}=1$ for any $y\in \cal X$;

(3) If $B_1 $   is a linear operator and $a>0$ is a real number, then $ {\|aB_1\|_y}\leq {\|B_1\|_{(ay)}} {\|aI\|_y}$;

(4) If $ {F}(B_1(y))\geq  {F}(B_2(y))$ for any $y\in \cal X$, then
$$ {\|B_1\|}_y\geq  {\|B_2\|}_y;$$

(5) If $ {F}(B_1(y)+ B_2(y))= {F}(B_1(y))+  {F}(B_2(y))$ for any $y\in \cal X$, then
$$ {\|B_1+ B_2\|_y}\leq {\|B_1\|}_y+  {\|B_1\|}_y;$$

(6) If $ {F}(B_1(y)- B_2(y))= {F}(B_1(y))-  {F}(B_2(y))\geq 0$ for any $y\in \cal X$, then
$$0\leq  {\|B_1\|}_y-  {\|B_2\|}_y\leq  {\|B_1- B_2\|_y};$$

(7) If  $(B_1B_2)(y)=B_1(B_2(y))$ for any $y\in \cal X$, and there exist $ {\|B_1\|_{ B_2(y)}}$ and $ {\|B_1\|_y}$, then
$$ {\|B_1 B_2\|_y}\leq {\|B_1\|_{B_2(y)}}\|B_2\|_y;$$

(8) If  ${\cal F(X)}_1\subseteq {\cal F(X)}_2\subseteq \cal P(X)$,
\begin{equation*}
 {\|B\|_y}_{{,\cal F(X)}_1}= \sup_{ {F}\in {\cal F(X)}_1}\frac{ {F}(B(y))}{ {F}(y)}
\end{equation*}
and
\begin{equation*}
 {\|B\|_y}_{{,\cal F(X)}_2}= \sup_{ {F}\in {\cal F(X)}_2}\frac{ {F}(B(y))}{ {F}(y)},
\end{equation*}
where  ${\|B\|_y}_{{,\cal F(X)}_1}$ represents the fuzzy rate of  the operator $B$ at the  point $y\in \cal X$ on ${\cal F(X)}_1$ and ${\|B\|_y}_{{,\cal F(X)}_2}$ on ${\cal F(X)}_2$,  then for any $y\in \cal X$,
$$ {\|B\|_y}_{{,\cal F(X)}_1}\leq  {\|B\|_y}_{{,\cal F(X)}_2}.$$
\end{thm}

\textbf{Proof.} (1) It follows that $ {\|B\|_y}\geq 0$ for any $B\in \cal B(X)$ from Definition 3.1. On the other hand, if $ {\|B\|_y}= 0$, then
for any $ {F}\in \cal F(X)$, $ {F}(B(y))=0$. It is false because there exists a membership function  $ {F}\in \cal F(X)$ where $ {F}(z)=0.5$ when $z=B(y)$ and $ {F}(z)=0$ when $z\neq B(y)$.

(2) Because  $I$ is an identity operator, we obtain
$$
 {\|I\|_y}= \sup_{ {F}\in \cal F (X)}\frac{ {F}(I(y))}{ {F}(y)}=\sup_{ {F}\in \cal F (X)}\frac{ {F}(I(y))}{ {F}(y)}=1
$$
 for all  $y\in \cal X$.

(3) If $B_1 $ is a linear operator and $a>0$ is a real number, then
\begin{eqnarray*}
 {\|aB_1\|_y}&=&\sup_{ {F}\in \cal F (X)}\frac{ {F}(aB_1(y))}{ {F}(y)}=\sup_{ {F}\in \cal F (X)}(\frac{ {F}(B_1(ay))}{ {F}(ay)}\frac{ {F}(ay)}{ {F}(y)})\nonumber\\
&\leq&\sup_{ {F}\in \cal F (X)}\frac{ {F}(B_1(ay))}{ {F}(ay)}\sup_{ {F}\in \cal F (X)}\frac{ {F}(ay)}{ {F}(y)}
= {\|B_1\|_{(ay)}} {\|aI\|_y}.
\end{eqnarray*}

(4) If $ {F}(B_1(y))\geq  {F}(B_2(y))$ for any $y\in \cal X$, then
\begin{eqnarray*}
 {\|B_1\|_y}=\sup_{ {F}\in \cal F (X)}\frac{ {F}(B_1(y))}{ {F}(y)}
\geq\sup_{ {F}\in \cal F (X)}\frac{ {F}(B_2(y))}{ {F}(y)}
= {\|B_2\|_y}.
\end{eqnarray*}

(5) If $ {F}(B_1(y)+ B_2(y))= {F}(B_1(y))+  {F}(B_2(y))$ for any $y\in \cal X$, we get
\begin{eqnarray*}
 {\|B_1+ B_2\|_y}&=&\sup_{ {F}\in \cal F (X)}\frac{ {F}((B_1+ B_2)(y))}{ {F}(y)}\nonumber\\
&\leq&\sup_{ {F}\in \cal F (X)}\frac{ {F}(B_1(y))}{ {F}(y)}+\sup_{ {F}\in \cal F (X)}\frac{ {F}(B_2(y))}{ {F}(y)}
= {\|B_1\|_y}+  {\|B_2\|_y}.\nonumber\\
\end{eqnarray*}

(6) Let $ {F}(B_1(y)- B_2(y))= {F}(B_1(y))-  {F}(B_2(y))\geq 0$ for any $y\in \cal X$, which means

 $$ \frac{ {F}(B_1(y))-  {F}(B_2(y))}{ {F}(y)}\geq 0$$
 and
\begin{eqnarray*}
 {\|B_1\|_y}&=&\sup_{ {F}\in \cal F (X)}\frac{ {F}(B_1(y))}{ {F}(y)}
=\sup_{ {F}\in \cal F (X)}\frac{[ {F}(B_1(y))- {F}(B_2(y))]+ {F}(B_2(y))}{ {F}(y)}\nonumber\\
&\leq &\sup_{ {F}\in \cal F (X)}\frac{ {F}(B_1(y))- {F}(B_2(y))}{ {F}(y)}+\sup_{ {F}\in \cal F (X)}\frac{ {F}(B_2(y))}{ {F}(y)}
= {\|B_1-B_2\|_y}+ {\|B_2\|_y}.
\end{eqnarray*}
Therefore, $0\leq  {\|B_1\|_y}-  {\|B_1\|_y}\leq  {\|B_1-B_2\|_y}$ holds by (4).

(7) Set $(B_1B_2)(y)=B_1(B_2(y))$ for any $y\in \cal X$. Then there exist $ {\|B_1\|_{ B_2(y)}}$ and $ {\|B_1\|_y}$ such that
\begin{eqnarray*}
 {\|B_1B_2\|_y}&=&\sup_{ {F}\in \cal F (X)}\frac{ {F}(B_1(B_2(y)))}{ {F}(y)}
=\sup_{ {F}\in \cal F (X)}\frac{ {F}(B_1(B_2(y)))}{ {F}(B_2(y))}\frac{ {F}(B_2(y))}{ {F}(y)}\nonumber\\
&\leq &\sup_{ {F}\in \cal F (X)}\frac{ {F}(B_1(B_2(y)))}{ {F}(B_2(y))}\sup_{ {F}\in \cal F (X)}\frac{ {F}(B_2(y))}{ {F}(y)}
= {\|B_1\|_{ B_2(y)}} {\|B_1\|_y}.
\end{eqnarray*}

(8) It is clear that the result holds by  ${\cal F(X)}_1\subseteq {\cal F(X)}_2\subseteq \cal P(X)$ and Definition 3.1.

The following is also a property of the  fuzzy rate of  the operator $B$ based on its basic properties.

\begin{cor}
Let $\cal X$ be a linear space, $B:\cal X \to \cal X $  be  an operator, $\cal {B}( X)$$=\{B|B: \cal X\to \cal X\}$, $F:\cal X$ $ \to [0, 1]$ be a membership function over $\cal X$, $\cal P(X)$$=\{F|F:\cal X$$ \to [0, 1]\}$  be a collection of all membership functions over $\cal X$. For any $\emptyset\neq \cal F (X)\subseteq \cal P(X)$, if $B^n(y)=B^{n-1}(B(y))$ for $n=1, 2, \cdot\cdot\cdot$, there exist $ {\|B\|}_{B^{k-1}(y)}$ for $k=1, 2, \cdot\cdot\cdot, n$, such that
\begin{equation}
\frac{1}{\prod_{1\leq k\leq n} {\|B\|}_{B^{k-1}(y)}}\leq\frac{1}{ {\|B^n\|}_y},
\end{equation}
where $B^0=I$ is an identity operator.
\end{cor}

\textbf{Proof.}  We have
 \begin{eqnarray*}
 {\|B^n\|_y}&=&\sup_{ {F}\in \cal F (X)}\frac{ {F}(B(B^{n-1}(y)))}{ {F}(y)}
=\sup_{ {F}\in \cal F(X)}\frac{ {F}((B(B^{n-1}(y)))}{ {F}(B^{n-1}(y))}\frac{ {F}(B^{n-1}(y))}{ {F}(y)}\nonumber\\
&\leq &\sup_{ {F}\in \cal F(X)}\frac{ {F}((B(B^{n-1}(y))))}{ {F}(B^{n-1}(y))}\sup_{ {F}\in \cal F (X)}\frac{ {F}(B^{n-1}(y))}{ {F}(y)}\nonumber\\
&=& {\|B\|}_{B^{n-1}(y)}\sup_{ {F}\in \cal F (X)}\frac{ {F}(B^{n-1}(y))}{ {F}(y)}
\leq \cdot\cdot\cdot
\leq \prod_{1\leq k\leq n} {\|B\|}_{B^{k-1}(y)}.
\end{eqnarray*}
It follows that the result (3.3) holds.

In what follows, we will apply the above properties to prove a new Fixed Point Existence Theorem.

 \section{ Applications--a new Fixed Point Existence Theorem}

Fixed Point theory is very important and most generally useful  one in classical function analysis.
In this section, we prove a new Fixed Point Existence Theorem with the fuzzy rate of  the operator $B$ as its application. First, we have

 \begin{lem}
Let $\cal X$ be a linear space, $B:\cal X \to \cal X $  be  an operator, $\cal {B}( X)$$=\{B|B: \cal X\to \cal X\}$, $F:\cal X$ $ \to [0, 1]$ be a membership function over $\cal X$, $\cal P(X)$$=\{F|F:\cal X$$ \to [0, 1]\}$  be a collection of all membership functions over $\cal X$. Let

\begin{equation*}
 {\|B\|_y}= \sup_{ {F}\in \cal F(X)}\frac{ {F}(B(y))}{ {F}(y)}<+\infty
\end{equation*}
for $\emptyset\neq \cal F (X)\subseteq \cal P(X)$. If for $\delta\in (0, 1]$, there exists a natural number $N$ such that $ {\|B^n\|_y}\geq \delta$ as $n\geq N$, then there exists a $ {F}_0\in \cal F (X)$ such that
\begin{equation}
 {F}_0(B^n(y))= {F}_0(B^{n-1}(y)),
\end{equation}
or
\begin{equation}
 {F}_0(B(B^{n-1}(y)))= {F}_0(B^{n-1}(y)),
\end{equation}
for $n\geq N$.
\end{lem}
\par

\textbf{Proof.} Note that  $B^n(y)=B^{n-1}(B(y))$ for $n=1, 2, \cdot\cdot\cdot$.

  If  $\delta\in (0, 1]$, there exists a  natural number $N$ such that $ {\|B^n\|_y}\geq \delta$ as $n\geq N$. Then we know
\begin{equation*}
1\leq \frac{1}{\prod_{1\leq k\leq n} {\|B\|}_{B^{k-1}(y)}}\leq\frac{1}{ {\|B^n\|}_y}\leq \frac{1}{\delta}<+\infty
\end{equation*}
for $n\geq N$ by (2.3), and
\begin{equation*}
0\leq \sum^n_{k=1}-\ln {\|B\|}_{B^{k-1}(y)}\leq -\ln {\|B^n\|}_y\leq -\ln\delta<+\infty,
\end{equation*}
hence $\lim_{k\rightarrow +\infty}\ln {\|B\|}_{B^{k-1}(y)}=0$ and $\lim_{k\rightarrow +\infty} {\|B\|}_{B^{k-1}(y)}=1$\cite{Zorich25}.

It follows that for any natural number $m$ there exists a $M$, as $k> max\{M, N\}$ such that
\begin{equation*}
1-\frac{1}{m}<  {\|B\|}_{B^{k-1}(y)}< 1+\frac{1}{m}.
\end{equation*}

Since
\begin{equation*}
 {\|B\|}_{B^{k-1}(y)}= \sup_{ {F}\in \cal F(X)}\frac{ {F}(B(B^{k-1}(y)))}{ {F}(B^{k-1}(y))}< +\infty,
\end{equation*}
 there exists a membership function  $ {F}_0\in \cal F(X)$ such that
\begin{equation*}
1-\frac{1}{m}-\frac{1}{m}<  {\|B\|}_{B^{k-1}(y)}-\frac{1}{m}<\frac{ {F}_0(B(B^{k-1}(y)))}{ {F}_0(B^{k-1}(y))}<  {\|B\|}_{B^{k-1}(y)}+\frac{1}{m}< 1+\frac{1}{m}+\frac{1}{m}.
\end{equation*}

Letting  $m\rightarrow +\infty$, we obtain

\begin{equation*}
\frac{ {F}_0(B(B^{k-1}(y)))}{ {F}_0(B^{k-1}(y))}=1,
\end{equation*}
 that's to say, $ {F}_0(B(B^{k-1}(y)))= {F}_0(B^{k-1}(y))$.

Then, we give the definition of a quasi-fixed point of the operator $B$ with respect to the membership function  $ {F}$.

 \begin{defn}
Let $\cal X$ be a linear space, $B:\cal X \to \cal X $  be  an operator, $\cal {B}( X)$$=\{B|B: \cal X\to \cal X\}$, $F:\cal X$ $ \to [0, 1]$ be a membership function over $\cal X$, $\cal P(X)$$=\{F|F:\cal X$$ \to [0, 1]\}$  be a collection of all membership functions over $\cal X$.
For $y\in \cal X$, if there exists  $ {F}\in \cal P(X)$ such that $ {F}(B(y))= {F}(y)$, then  $y$ is called a quasi-fixed point of the operator $B$ with respect to $ {F}$.
\end{defn}

By  the proof of Lemma 4.1,  $B^{k-1}(y)$ is  a quasi-fixed point of the operator $B$ with respect to $ {F}$.

Now, Fixed Point Existence Theorem with respect to $ {F}$ is presented.

 \begin{thm}(new Fixed Point Existence Theorem)
Let $\cal X$ be a linear space, $B:\cal X \to \cal X $  be  an operator, $\cal {B}( X)$$=\{B|B: \cal X\to \cal X\}$, $F:\cal X$ $ \to [0, 1]$ be a membership function over $\cal X$, $\cal P(X)$$=\{F|F:\cal X$$ \to [0, 1]\}$  be a collection of all membership functions over $\cal X$.
Presume

 $$ {\|B\|_y}= \sup_{ {F}\in \cal F(X)}\frac{ {F}(B(y))}{ {F}(y)}<+\infty$$ for $\emptyset\neq \cal F (X)\subseteq \cal P(X).$ If  $\delta\in (0, 1],$ there exists a natural number $N$ such that $ {\|B^n\|_y}\geq \delta$ as $n\geq N$. Then if there exists an injection functional $ {F}_0\in \cal F (X)$ such that
\begin{equation}
 {F}_0(B^n(y))= {F}_0(B^{n-1}(y)),
\end{equation}
 $B^{n-1}(y)$ is a fixed point of the operator $B$ with respect to $ {F}$, and  $y$ is a fixed point of the operator $B^n$ with respect to $ {F}$ for $n\geq N$.
\end{thm}
\par
\textbf{Proof.} It follows directly that the result holds from (4.2) and the injective condition of functional $ {F}_0\in \cal F(X)$.

Like the classical fixed point theory applied to differential equations, we believe that the Fixed Point Existence Theorem with respect to fuzzy set $ {F}$ might be applied to fuzzy equations or  fuzzy differential equations. They are worth further studying in the future.

\section{Conclusions}

In this work, we have obtained the following results:

$\bullet$ The fuzzy rate of an operator in linear spaces is introduced
and some properties and basic principles of the fuzzy rate are  studied.

$\bullet$ The fuzzy rate of an diagonal matrix $B$ in a plane is discussed.

$\bullet$ A  new fixed point existence theorem is proved.

\section{Acknowledgments}

$\bullet$  This work was supported by Natural Science Foundation Project of Chongqing(Grant No. cstc2019jcyj-msxmX0716).

$\bullet$ \textbf{Competing interests}  The authors declare that they have no competing interests regarding the publication of this article.

\noindent
\\

\end{document}